\documentclass[11pt,reqno]{amsart}
\usepackage[usenames]{color}
\usepackage{amsmath,verbatim,graphicx,epstopdf,enumerate}
\newcommand{\Beq}{\begin{equation}}
\newcommand{\Eeq}{\end{equation}}
\newcommand{\beq}{\begin{equation*}}
\newcommand{\eeq}{\end{equation*}}
\newcommand{\bal}{\begin{align}}
\newcommand{\eal}{\end{align}}

\newcommand{\bp}{\begin{prob}}
	\newcommand{\ep}{\end{prob}}
\newcommand{\bpr}{\begin{proof}}
	\newcommand{\epr}{\end{proof}}

\newcommand{\bel}[1]{\begin{equation}\label{#1}}
\newcommand{\ee}{\end{equation}}

\newcommand{\rr}{\mathbb{R}}

\newtheorem{theorem}{Theorem}[section]

\newtheorem{lemma}[theorem]{Lemma}
\newtheorem{prop}[theorem]{Proposition}

\theoremstyle{definition}

\newtheorem{remark}[theorem]{Remark}

\title[An Inverse problem for the wave operator]{An inverse problem for the  wave equation with source and receiver at distinct points}
\author[Vashisth]{Manmohan Vashisth}
\address{Beijing Computational Science Research Center, Beijing 100193, China.
	\newline
	\indent E-mail: {\tt manmohanvashisth@gmail.com}}

\begin{document}
	\maketitle
	\begin{abstract}
		We consider the inverse problem of determining the density coefficient appearing in the wave equation from \emph{separated} point source and point receiver data. Under some assumptions on the coefficients, we
		prove uniqueness results.
	\end{abstract}
	
	

	\section{Introduction}
	We address the inverse problem of determining the density coefficient of a medium by probing it with an external point source and by measuring the responses at a single point for a certain period of time. 
	
	More precisely, consider the following initial value problem (IVP), where $\Box =\partial_{t}^{2}-\Delta_{x}$ denotes the wave operator:
	
	\begin{equation}\label{wave equation with potential}
	\begin{aligned}
	(\Box-q(x))u(x,t)&=\delta(x,t),\ \ (x,t)\in\rr^{3}\times\rr\\
	u(x,t)|_{t<0}&=0,\ \ \ \quad \quad \quad \  x\in\rr^{3}.
	\end{aligned}
	\end{equation}

	In  Equation \eqref{wave equation with potential}, we  assume that the coefficient $q$ is real-valued and is a   {\color{black}$C^{3}(\rr^{3})$} function. The  inverse problem we address  is the unique determination of the coefficient $q$ from the knowledge of $u(e,t)$ where $e=(1,0,0)$ for $t\in[0,T]$ with $T>1$. Motivation for studying such problems  arises in geophysics see \cite{Symes} and references therein. Geophysicist determine properties of the earth structure by sending waves from the surface of the earth and measuring the corresponding scattered responses. Note that in the problem we consider here, the point source is located at the origin, whereas the responses are measured at a different point.
	Since the  given data depends on one variable whereas the coefficient to be determined depends on three variables,  some additional restrictions on the coefficient $q$ are required to make the inverse problems tractable. 
	
	
	{\color{black}There are several results related to point source inverse problems involving the wave equation. We briefly mention here the details of works which are closely related to the problem studied in this article. Romanov in \cite{Romanov} considered the problem of determining the damping and density coefficients which are constant outside a bounded, simply connected domain $D\subset\mathbb{R}^{3}$. By using the expression for fundamental solution, he reduced the problem to an integral geometry problem (whose solution was known by \cite{LDS}), which gives the determination of these coefficients in  $D$ when source and receiver are moving in a plane ( say $M$) chosen in such a way that $(i)$ $M\cap\overline{D}=\emptyset$ and $(ii)$ the line segment joined by source and receiver lies completely outside $\overline{D}$. Rakesh in \cite{Rakesh single-coincident source-receiver pair} studied the problem of determining $q$ from the knowledge  of $u(0,t)$ for $t\in [0,T]$ and he proved the uniqueness  for coefficients which  
		are either comparable or  radially symmetric with respect to a point different from the source location. The above mentioned works are related to point source hyperbolic inverse problems with under-determined data. We also mention some related works for the point source hyperbolic inverse problems with formally determined or with over determined data. In \cite{Rakesh inversion of spherically symmetric potential} Rakesh proved the unique determination of the radially symmetric coefficient $q(x)$   appearing in \eqref{wave equation with potential} when  $u(a,t)$ is known for all $a\in \mathbb{S}^{2}$ and $t\in [0,T]$.
		Rakesh and Sacks in \cite{Rakesh and Sacks uniqueness angular controlled potential}  established the  uniqueness for angular controlled coefficient appearing in \eqref{wave equation with potential} from the knowledge of $u(a,t)$ and $u_{r}(a,t)$ for all $a\in \mathbb{S}^{2}$ and $t\in [0,T]$ where $u_{r}$ denotes the derivative with respect to $r=\lvert x\rvert$. The problem considered in  \cite{Rakesh and Sacks uniqueness angular controlled potential} can be seen as an extension of the work \cite{Rakesh inversion of spherically symmetric potential} to a set of more general coefficients which is strictly bigger than the set of radial functions but \cite{Rakesh and Sacks uniqueness angular controlled potential} requires more information than \cite{Rakesh inversion of spherically symmetric potential}. In \cite{Rakesh and Uhlmann inverse backscattering} the problem of  determining the density coefficient $q$ with angular controlled is studied. They proved the uniqueness of these coefficients from the knowledge of $u^{a}(a,t)$ for all $a\in \mathbb{S}^{2}$ and $t\in [0,T]$, where $u^{a}(x,t)$ denote the  solution to \eqref{wave equation with potential} when source is located at $a\in\mathbb{S}^{2}$. For more works related to the problem studied in this article, we refer to \cite{Burridge integral equation,Sacks and Symes uniqueness and continuous dependence,Rakesh and Symes uniqueness result,Stefanov uniqueness,Romanov,Rakesh impedance inversion,Rakesh layered medium point source,Klibanov partial point source,Li partial point source,Romanov book integral geometry} and references therein.} For
	the one dimensional inverse problems related to the problem studied here, we refer to \cite{Bube and Burridge 1-D,Bukhgem and Klibanov uniqueness,Rakesh and Sacks impedance inversion, Rakesh Webster's equation}.
	

	
	
	We now state the main results of this article.
	
	\begin{theorem}\label{Separated Result 1}
		Suppose {\color{black}$q_{i}\in C^{3}(\rr^{3})$}, $i=1,2$  with $q_{1}(x)\geq q_{2}(x)$ {\color{black}{for all}} $x\in\mathbb{R}^{3}$. Let $u_{i}(x,t)$ be the solution to the IVP
		\begin{align*}
		\begin{aligned}
		(\Box-q_{i}(x))u_{i}(x,t)&=\delta(x,t),\ \ \ \ (x,t)\in\rr^{3}\times\rr\\
		u_{i}(x,t)|_{t<0}&=0, \ \ \ \ \ \ \ \ \ \ \ \ \ \   x\in\rr^{3}{\color{black}.}
		\end{aligned}
		\end{align*} 
		If $u_{1}(e,t)=u_{2}(e,t)$, for all $t\in[0,T]$ {\color{black} where} $T>1$ and $e=(1,0,0)$, then  $q_{1}(x)=q_{2}(x)$ for all $x$ with $|e-x|+|x|\leq T$.
	\end{theorem}
	\begin{theorem}\label{Separated Result 2}
		Suppose {\color{black}$q_{i}\in C^{3}(\rr^{3})$}, $i=1,2$ with $q_{i}(x)=a_{i}(|x|+|x-e|)$ with $e=(1,0,0)$, for some {\color{black}$C^{3}$} functions $a_{i}$ on $(1-\epsilon,\infty)$ for some $0<\epsilon<1$.  Let $u_{i}$ be the solution to the IVP
		\begin{align*}
		\begin{aligned}
		(\Box-q_{i}(x))u_{i}(x,t)&=\delta(x,t),\ \ \ \ (x,t)\in\rr^{3}\times\rr\\
		u_{i}(x,t)|_{t<0}&=0, \ \ \ \ \ \ \ \ \ \ \ \ \ \ \ \  x\in\rr^{3} {\color{black}.}
		\end{aligned}
		\end{align*} 
		If $u_{1}(e,t)=u_{2}(e,t)$, for all $t\in[0,T]$ {\color{black} where} $T>1$ and $e=(1,0,0)$, then $q_{1}(x)=q_{2}(x)$ for all $x$ with $|e-x|+|x|\leq T$.
	\end{theorem}
	{\begin{remark}\label{remark for T larger than 1}
			From Proposition $2.1$, the solution $u(x,t)$ of $(1.1)$ is supported in $t\geq \lvert x\rvert$  hence $u(e,t)=0$ for $t<1$. So $u(e,t)$ has no information about $q$ if $t<1$, hence we require $T>1$ in Theorems \ref{Separated Result 1} and \ref{Separated Result 2}. Further note that ellipsoids $\lvert x\rvert + \lvert x-e\rvert \leq t$ are empty if $t<1$.
	\end{remark}}

	To the best of our knowledge, our  results, Theorem{\color{black}s} \ref{Separated Result 1} and  \ref{Separated Result 2},  which treat \emph{separated} source and  receiver, have not been studied earlier. Our result generalize the work \cite{Rakesh single-coincident source-receiver pair}, who considered the aforementioned inverse problem but with \emph{coincident} source and receiver;  see also \cite{Stefanov uniqueness}.

	The proofs of the above theorems are based on an integral identity derived using the solution to an adjoint problem as used in \cite{Santosa Syemes} and \cite{Stefanov uniqueness}.  Recently this idea was used in \cite{Rakesh and Uhlmann inverse backscattering} as well.

	The article is organized as follows. In Section \ref{Prelim}, we state the existence and uniqueness results for the solution of  Equation \eqref{wave equation with potential}, the proof of which is given in \cite{Friedlander book,Levrentev Romanov and Shishatski book ill-posed,Romanov book integral geometry}. Section \ref{Proof Th 1 and Th 2} contains the proofs of Theorems \ref{Separated Result 1} and \ref{Separated Result 2}. 
	\section{Preliminaries}\label{Prelim}
	\begin{prop}\cite[pp.139,140]{Friedlander book}\label{Damping Prop}
		Suppose {\color{black}$q$ is a $C^{3}$} function on $\rr^{3}$
		and $u(x,t)$ satisfies the following IVP
		\begin{equation}\label{Hyperbolic PDE}
		\begin{aligned}
		(\Box -q(x))u(x,t)&= \delta(x,t),  \ (x,t)\in\rr^{3}\times\rr\\
		u(x,t)|_{t<0}\ &=\ 0, \quad \quad \quad \ \ \  x\in \rr^{3}
		\end{aligned}
		\end{equation}
		then $u(x,t)$ is given by 
		\begin{equation}\label{Fundamental solution to hyperbolic equation}
		u(x,t) =\frac{\delta(t-\lvert x\rvert)}{4\pi\lvert x\rvert} + v(x,t)
		\end{equation}
		where $v(x,t)=0$ for $t<|x|$ and in the region $t>\lvert x\rvert$, $v(x,t)$ is a $C^{2}$
		solution of the characteristic boundary value problem (Goursat Problem)
		\begin{equation}\label{Goursat problem for zeroth order}
		\begin{aligned}
		(\Box - q(x))v(x,t)&=0,\ \ \  \ t>\lvert x\rvert\\
		v(x,\lvert x\rvert)&=\frac{1}{8\pi}\int\limits_{0}^{1}q(sx)ds.
		\end{aligned}
		\end{equation}
		
	\end{prop}
	We will use the following version of this proposition. Consider the following IVP
	\begin{align}\label{IVP with source at e}
	\begin{aligned}
	(\Box -q(x))U(x,t)&= \delta(x-e,t),  \ (x,t)\in\rr^{3}\times\rr\\
	U(x,t)|_{t<0}\ &=\ 0, \quad \qquad \qquad \ \  x\in \rr^{3}.
	\end{aligned}
	\end{align}
	Now we have 
	\begin{align}\label{expression for solution to IVP with source at e}
	\begin{aligned}
	U(x,t)=\frac{\delta(t-\lvert x-e\rvert)}{4\pi\lvert x-e\rvert}+V(x,t)
	\end{aligned}
	\end{align}
	where $V(x,t)=0$ for $t<\lvert x-e\rvert $ and for $t>\lvert x-e\rvert$, $V(x,t)$ is a $C^{2}$ solution to the following Goursat Problem 
	\begin{align}\label{equation for V when source is at e}
	\begin{aligned}
	(\Box - q(x))V(x,t)&=0,\ \ \  \ t>\lvert x-e\rvert \\
	V(x,\lvert x-e\rvert)&=\frac{1}{8\pi}\int\limits_{0}^{1}q(sx+(1-s)e)ds.
	\end{aligned}
	\end{align}
	We can see this by translating source  by $-e$ in Equation \eqref{IVP with source at e} and using the above proposition.
	
	\renewcommand{\thesection}{\large 3}
	\section{\large  Proof of Theorems \ref{Separated Result 1} and \ref{Separated Result 2}}\label{Proof Th 1 and Th 2}
	In this section, we prove {\color{black}Theorems} \ref{Separated Result 1} and \ref{Separated Result 2}. We will first show the following three lemmas which will be used in the proof of the main results.
	\begin{lemma}\label{Separated Integral Identity}
		Suppose $q_{i}'s$ $i=1,2$ be {\color{black}$C^{3}$} real-valued functions on $\rr^{3}$.
		Let $u_{i}$ be the solution to Equation \eqref{wave equation with potential} with $q=q_{i}$ and denote $u(x,t):=u_{1}(x,t)-u_{2}(x,t)$ and 
		$q(x):= q_{1}(x)-q_{2}(x)$. Then we have the following  integral identity  
		\begin{equation}\label{integral Identity}
		u(e,\tau)=\int\limits_{\rr}\int\limits_{\rr^{3}}q(x)u_{2}(x,t)U(x,\tau-t)dxdt, \mbox{ for all } \tau\in\rr
		\end{equation}
		where $U(x,t)$ is the solution to the following IVP
		\begin{equation}\label{equation for U}
		\begin{aligned}
		(\Box-q_{1}(x))U(x,t)&=\delta(x-e,t),\ \ (x,t)\in\rr^{3}\times\rr\\
		U(x,t)|_{t<0}&=0,\ \ \ \quad \quad \quad\quad\quad \  x\in\rr^{3}.
		\end{aligned}
		\end{equation}
		\begin{proof}
			Since each $u_{i}$ for $i=1,2$ satisfies the following IVP,
			\begin{align*}
			\begin{aligned}
			(\Box - q_{i}&(x))u_{i}(x,t) = \delta(x,t),  \ (x,t)\in\rr^{3}\times\rr\\
			&u_{i}(x,t)|_{t<0}\ =\ 0, \quad \quad \ \quad\quad  x\in \rr^{3},
			\end{aligned}
			\end{align*}
			we have that $u$ satisfies the following IVP
			\begin{equation}\label{equation for u}
			\begin{aligned}
			(\Box - q_{1}&(x))u(x,t) = q(x)u_{2}(x,t),  \ (x,t)\in\rr^{3}\times\rr\\
			&u(x,t)|_{t<0}\ =\ 0, \quad \quad \  x\in \rr^{3}.
			\end{aligned}
			\end{equation}
			
			Now since 
			\begin{align*}
			\begin{aligned}
			u(e,\tau)=\int\limits_{\rr^{3}}\int\limits_{\rr}u(x,t)\delta(x-e,\tau-t)dtdx,
			\end{aligned}
			\end{align*}
			using $\eqref{equation for U}$,  we have
			\begin{align*}
			\begin{aligned}
			u(e,\tau)=\int\limits_{\rr^{3}}\int\limits_{\mathbb{\rr}}u(x,t)(\Box -q_{1}(x))U(x,\tau-t)dtdx.
			\end{aligned}
			\end{align*}
			
			Now by using integration by parts and  Equations \eqref{equation for U} and \eqref{equation for u}, also taking into account that $u(x,t)=0$ for $t<\lvert x\rvert $ and that $U(x,t)=0$ for $\lvert x-e\rvert>t$,
			we get 
			\begin{align*}
			\begin{aligned}
			u(e,\tau)&=\int\limits_{\rr^{3}}\int\limits_{\rr}U(x,\tau-t)(\Box-q_{1}(x))u(x,t)dtdx\\
			&=\int\limits_{\rr^{3}}\int\limits_{\rr}q(x)u_{2}(x,t)U(x,\tau-t)dtdx.
			\end{aligned}
			\end{align*}
			This completes the proof of the lemma.
		\end{proof}
	\end{lemma}
	
	\begin{lemma}
		Suppose $q_{i}'s$ are as in Lemma \ref{Separated Integral Identity} and $u_{i}$ is the solution to Equation \eqref{wave equation with potential} with $q=q_{i}$ and if $u(e,t):=(u_{1}-u_{2})(e,t)=0$ for all $t\in[0,T]$, then there exists a constant $K>0$ depending on the bounds on $v_{2}$, $V$and $T$ such that the following inequality holds
		\begin{equation}\label{Estimate}
		\Bigg \lvert \int\limits_{\lvert x-e\rvert+\lvert x\rvert=2\tau}\frac{q(x)}{\lvert 2\tau x-\lvert x\rvert e\rvert}dS_{x}\Bigg\rvert \leq K \int\limits_{\lvert x-e\rvert+\lvert x\rvert\leq 2\tau}\frac{\lvert q(x)\rvert}{\lvert x\rvert \lvert x-e\rvert}dx,\ \forall\tau\in(1/2,T/2].
		\end{equation}
		Here $dS_x$ is the surface measure on the ellipsoid $\lvert x-e\rvert+\lvert x\rvert=2\tau$ and  $v_{2}$, $V$  are solutions to the Goursat problem (see Equations \eqref{Goursat problem for zeroth order} and \eqref{equation for V when source is at e}) corresponding to $q=q_{i}$.
		\begin{proof}
			From Lemma \ref{Separated Integral Identity}, we have 
			\begin{align*}
			u(e,2\tau)=\int\limits_{\rr}\int\limits_{\rr^{3}}q(x)u_{2}(x,t)U(x,2\tau-t)dxdt, \mbox{ for all } \tau\in\rr.
			\end{align*}
			Now since $u(e,2\tau)=0$ for all $\tau\in[0,T/2]$, and using Equations \eqref{Fundamental solution to hyperbolic equation} and \eqref{expression for solution to IVP with source at e}, we get
			\begin{align*}
			\begin{aligned}
			0&=\int\limits_{\rr^{3}}\int\limits_{\rr}q(x)\frac{\delta(t-\lvert x\rvert)\delta(2\tau-t-\lvert x-e\rvert)}{16\pi^{2}\lvert x\rvert\lvert x-e\rvert}dtdx\\
			&\ \ + \int\limits_{\rr^{3}}\int\limits_{\rr}q(x)\frac{\delta(t-\lvert x\rvert)V(x,2\tau-t)}{4\pi\lvert x\rvert}dtdx\\
			&\ \ +
			\int\limits_{\rr^{3}}\int\limits_{\rr}q(x)\frac{\delta(2\tau-t-\lvert x-e\rvert)}{4\pi \lvert x-e\rvert}v_{2}(x,t)dtdx\\ &\ \ +\int\limits_{\rr^{3}}\int\limits_{\rr}q(x) V(x,2\tau-t)v_{2}(x,t)dtdx. 
			\end{aligned}
			\end{align*}
			
			Now using the fact that $v_{2}(x,t)=0$ for $t<\lvert x\rvert$, $V(x,t)=0$ for $t<\lvert x-e\rvert$ and 
			\begin{align*}
			\int\limits_{\rr^{n}}\phi(x)\delta(P)dx=\int\limits_{P(x)=0}\frac{\phi(x)}{\lvert\nabla_{x}P(x)\rvert}dS_{x}
			\end{align*}
			where $dS_{x}$ is the surface measure on the surface $P=0$, we have that
			\begin{align*}
			\begin{aligned}
			0=&\frac{1}{16\pi^{2}}\int\limits_{\lvert x-e\rvert+\lvert x\rvert=2\tau}\frac{q(x)}{\lvert x\rvert \lvert x-e\rvert\lvert\nabla_{x}(2\tau-\lvert x\rvert-\lvert x-e\rvert)\rvert}dS_{x}\\
			& + \frac{1}{4\pi}\int\limits_{\lvert x\rvert+\lvert x-e\rvert\leq 2\tau}\frac{q(x)V(x,2\tau-\lvert x\rvert)}{\lvert x\rvert}dx\\
			&+\frac{1}{4\pi}\int\limits_{\lvert x\rvert+\lvert x-e\rvert \leq 2\tau}\frac{q(x)v_{2}(x,2\tau-\lvert x-e\rvert)}{\lvert x-e\rvert}dx\\
			&+\int\limits_{\lvert x\rvert+\lvert x-e\rvert\leq 2\tau}\int\limits_{\lvert x\rvert}^{2\tau -\lvert x-e\rvert}q(x)V(x,2\tau-t)v_{2}(x,t)dtdx.
			\end{aligned}
			\end{align*}
			For simplicity, denote 
			\begin{align*}
			F(\tau,x)&:=\frac{1}{4\pi}\Big(\lvert x-e\rvert V(x,2\tau-\lvert x\rvert)+\lvert x\rvert v_{2}(x,2\tau-\lvert x-e\rvert)\\
			&\ \ \  \ \ \ +4\pi\lvert x\rvert\lvert x-e\rvert\int\limits_{\lvert x\rvert}^{2\tau-\lvert x-e\rvert}V(x,2\tau-t)v_{2}(x,t)dt\Big)
			\end{align*}
			and using
			\begin{align*}
			\left\lvert\nabla_{x}(2\tau-|x|-|x-e|)\right\rvert=\left\lvert\frac{x}{|x|}+\frac{x-e}{|x-e|}\right\rvert=\left\lvert\frac{|x-e|x+(x-e)|x|}{|x||x-e|}\right\rvert.
			\end{align*}
			We have
			\begin{align*}
			\begin{aligned}
			\frac{1}{16\pi^{2}}\int\limits_{|x-e|+|x|=2\tau}\frac{q(x)}{|2\tau x-|x|e| }dS_{x}=-\int\limits_{|x|+|x-e|\leq 2\tau}\frac{q(x)}{|x||x-e|}F(\tau,x)dx.
			\end{aligned}
			\end{align*}
			Note that $\tau\in[0,T/2]$ with $T<\infty$. Now using the boundedness of $v_{2}$ and $V$ on compact subsets, we have $|F(\tau,x)|\leq K$ on $|x|+|x-e|\leq T$.
			
			Therefore, finally we have 
			\begin{align*}
			\begin{aligned}
			\Bigg|\int\limits_{|x-e|+|x|=2\tau}\frac{q(x)}{|2\tau x-|x|e|}dS_{x}\Bigg|\leq K \int\limits_{|x-e|+|x|\leq 2\tau}\frac{|q(x)|}{|x||x-e|}dx,\ \forall\tau\in(1/2,T/2].
			\end{aligned}
			\end{align*}
			The lemma is proved.
		\end{proof}
	\end{lemma}

	\begin{lemma}
		Consider the solid ellipsoid $|e-x|+|x|\leq r$, where $e=(1,0,0)$ and $x=(x_{1},x_{2},x_{3})$, then we have its parametrization in  prolate-{\color{black}spheroidal} co-ordinates $(\rho,\theta,\phi)$ 
		given by 
		\begin{equation}\label{Prolate}
		\begin{aligned}
		x_{1}&=\frac{1}{2}+\frac{1}{2}\cosh\rho\cos\phi\\
		x_{2}&=\frac{1}{2}\sinh\rho\sin\theta\sin\phi\\
		x_{3}&=\frac{1}{2}\sinh\rho\cos\theta\sin\phi\\
		\end{aligned}
		\end{equation}
		with $\cosh\rho\leq r$, $\theta\in(0,2\pi)$, $\phi\in(0,\pi)$  and the surface measure $dS_{x}$ on $|e-x|+|x|=r$ 
		and volume  element $dx$ on $|e-x|+|x|\leq r$, are given by  
		\begin{equation}\label{Surface and Volume measure on ellipsoid}
		\begin{aligned}
		dS_{x}&=\frac{1}{4}\sinh\rho\sin\phi\sqrt{\cosh^{2}\rho -\cos^{2}\phi} d\theta d\phi, \ \\ & \mbox{with} \ \cosh\rho=r, \ \theta\in[0,2\pi] \ \mbox{and} \ \phi\in[0,\pi]\\
		dx&=\frac{1}{8} \sinh\rho\sin\phi(\cosh^{2}\rho-\cos^{2}\phi)d\rho d\theta d\phi, \\ &  \mbox{with} \ \cosh\rho\leq r, \ \theta\in[0,2\pi] \ \mbox{and} \ \phi\in[0,\pi]. 
		\end{aligned}
		\end{equation}
		\begin{proof}
			The above result is well known, but for completeness, we will give the proof. The solid ellipsoid $|e-x|+|x|\leq r$ in explicit form can be  written as 
			\begin{align*}
			\frac{(x_{1}-1/2)^{2}}{r^{2}/4}+\frac{x_{2}^{2}}{(r^{2}-1)/4}+\frac{x_{3}^{2}}{(r^{2}-1)/4}\leq 1.
			\end{align*}
			From this, we see that 
			\begin{align*}
			x_{1}&=\frac{1}{2}+\frac{1}{2}\cosh\rho\cos\phi\\
			x_{2}&=\frac{1}{2}\sinh\rho\sin\theta\sin\phi\\
			x_{3}&=\frac{1}{2}\sinh\rho\cos\theta\sin\phi\\
			\end{align*}
			with $\cosh\rho\leq r$, $\theta\in[0,2\pi]$ and $\phi\in[0,\pi]$. This proves the first part of the lemma.
			
			Now the parametrization of ellipsiod $|e-x|+|x|=r$, is given by 
			\begin{align*}
			F(\theta,\phi)=\left(\frac{1}{2}+\frac{1}{2}\cosh\rho\cos\phi,\frac{1}{2}\sinh\rho\sin\theta\sin\phi,\frac{1}{2}\sinh\rho\cos\theta\sin\phi\right)
			\end{align*}
			with  $\theta\in(0,2\pi)$, $\phi\in(0,\pi)$ and  $\cosh\rho=r$.
			
			Next, we have 
			\begin{align*}
			\begin{aligned}
			\frac{\partial F}{\partial\theta}&=\left(0,\frac{1}{2}\sinh\rho\cos\theta\sin\phi,-\frac{1}{2}\sinh\rho\sin\theta\sin\phi\right)\\
			\frac{\partial F}{\partial\phi}&=\left(-\frac{1}{2}\cosh\rho\sin\phi,\frac{1}{2}\sinh\rho\sin\theta\cos\phi,\frac{1}{2}\sinh\rho \cos\theta \cos\phi\right).
			\end{aligned}
			\end{align*}
			We have $dS_{x}=\left\lvert\frac{\partial F}{\partial\theta}\times \frac{\partial F}{\partial\phi}\right\rvert d\theta d\phi$,
			simple computation will gives us 
			\begin{align*}
			&dS_{x}=\frac{1}{4}\sinh\rho\sin\phi\sqrt{\cosh^{2}\rho -\cos^{2}\phi} d\theta d\phi, \\
			& \ \text{with} \ \cosh\rho=r, \ \theta\in(0,2\pi) \ and \ \phi\in(0,\pi).
			\end{align*}
			
			Last part of the lemma follows from change of variable formula, which is given by 
			\begin{align*}
			dx=\Big|\frac{\partial(x_{1},x_{2},x_{3})}{\partial(\rho,\theta,\phi)}\Big|d\theta d\phi d\rho; \ \text{with} \ \cosh\rho\leq r,\ \theta\in[0,2\pi]\ \text{and} \ \phi\in[0,\pi]
			\end{align*}
			where $\frac{\partial(x_{1},x_{2},x_{3})}{\partial(\rho,\theta,\phi)}$ is given by 
			\begin{align*}
			\frac{\partial(x_{1},x_{2},x_{3})}{\partial(\rho,\theta,\phi)}=\det
			\begin{vmatrix}
			\frac{\partial{x_{1}}}{\partial{\rho}} 
			& \frac{\partial{x_{1}}}{\partial{\theta}} 
			& \frac{\partial{x_{1}}}{\partial{\phi}}\\[3mm]
			\frac{\partial{x_{2}}}{\partial\rho} 
			& \frac{\partial{x_{2}}}{\partial{\theta}} 
			& \frac{\partial{x_{2}}}{\partial{\phi}} \\[3mm]
			\frac{\partial{x_{3}}}{\partial{\rho}}  
			& \frac{\partial{x_{3}}}{\partial{\theta}}  
			& \frac{\partial{x_{3}}}{\partial{\phi}}
			\end{vmatrix}.
			\end{align*}
			This gives 
			\begin{align*}
			&dx=\frac{1}{8} \sinh\rho\sin\phi(\cosh^{2}\rho-\cos^{2}\phi)d\rho d\theta d\phi, \\ 
			& \text{with} \ \cosh\rho\leq r, \ \theta\in[0,2\pi] \ \mbox{and} \ \phi\in[0,\pi].
			\end{align*} 
		\end{proof}
		
	\end{lemma}

	\renewcommand{\thesection}{\large 3}
	\subsection{Proof of Theorem \ref{Separated Result 1}}
	We first consider the surface integral in Equation $\eqref{Estimate}$ and denote it $Q(2\tau)$:
	\begin{equation}\label{Definition of Q}
	Q(2\tau):=\int\limits_{|x-e|+|x|=2\tau}\frac{q(x)}{|2\tau x-|x|e|}dS_{x}.
	\end{equation}
	We have
	\begin{align*}
	\begin{aligned}
	|2\tau x-|x|e|&=\lvert(2\tau x_{1}-|x|,2\tau x_{2},2\tau x_{3})\rvert=\sqrt{(2\tau x_{1}-|x|)^{2}+4\tau^{2}x_{2}^{2}+4\tau^{2}x_{3}^{2}}\\
	&=\sqrt{4\tau^{2}|x|^{2}+|x|^{2}-4\tau x_{1}|x|}{\color{black}.}
	\end{aligned}
	\end{align*}
	From Equation $\eqref{Prolate}$ and using the fact that $\cosh\rho=2\tau$, we have 
	\begin{align*}
	\begin{aligned}
	|2\tau x-|x|e|&=\frac{1}{2}\sqrt{(2\tau+\cos\phi)\{(4\tau^{2}+1)(2\tau+\cos\phi)-4\tau(1+2\tau\cos\phi)\}}\\
	&=\frac{1}{2}\sqrt{(2\tau+\cos\phi)(8\tau^{3}+4\tau^{2}\cos\phi+2\tau+\cos\phi-4\tau-8\tau^{2}\cos\phi)}\\
	&=\frac{1}{2}\sqrt{(2\tau+\cos\phi)(8\tau^{3}-4\tau^{2}\cos\phi-2\tau+\cos\phi)}\\
	&=\frac{1}{2}\sqrt{(4\tau^{2}-\cos^{2}\phi)(4\tau^{2}-1)}.
	\end{aligned}
	\end{align*}
	
	Using the above expression for $|2\tau x-|x|e|$ and
	Equation \eqref{Surface and Volume measure on ellipsoid}, we get
	\begin{align*}
	\begin{aligned}
	Q(2\tau)&=\frac{1}{8}\int\limits_{0}^{\pi}\int\limits_{0}^{2\pi}\frac{q(\rho,\theta,\phi)\sinh\rho \sin\phi\sqrt{\cosh^{2}\rho-\cos^{2}\phi}}{\sqrt{(4\tau^{2}-\cos^{2}\phi)(4\tau^{2}-1)}} d\theta d\phi,
	\end{aligned}
	\end{align*}
	where we have denoted 
	\[
	q(\rho,\theta,\phi)=q\left( \frac{1}{2}+\frac{1}{2}\cosh\rho \cos\phi,\frac{1}{2}\sinh\rho\sin\theta\sin\phi,\frac{1}{2}\sinh\rho \cos\theta \sin\phi\right).
	\]
	{\color{black}After} using $\cosh\rho=2\tau$, $\sinh\rho=\sqrt{4\tau^{2}-1}$ and $\rho = \ln(2\tau+\sqrt{4\tau^{2}-1})$, we get 
	\begin{equation}\label{expression for Q}
	\begin{aligned}
	Q(2\tau)=\int\limits_{0}^{\pi}\int\limits_{0}^{2\pi}q(\ln(2\tau+\sqrt{4\tau^{2}-1}),\theta,\phi)\sin\phi d\theta d\phi.
	\end{aligned}
	\end{equation}
	Now consider the integral 
	\[
	\int\limits_{|x-e|+|x|\leq 2\tau}\frac{q(x)}{|x||x-e|}dx.
	\]
	Again using Equations \eqref{Prolate} and \eqref{Surface and Volume measure on ellipsoid} in the above integral,  we have
	\begin{align*}
	\begin{aligned}
	\int\limits_{|x-e|+|x|\leq 2\tau}\frac{q(x)}{|x||x-e|}dx&= \frac{1}{2}\int\limits_{\cosh\rho\leq 2\tau}\int\limits_{0}^{\pi}\int\limits_{0}^{2\pi}q(\rho,\theta,\phi)\sinh\rho\sin\phi d\theta d\phi d\rho.\\
	\end{aligned}
	\end{align*}
	After substituting $\cosh\rho=r$ and $\rho=\ln(r+\sqrt{r^{2}-1})$,  we get 
	\begin{align*}
	\int\limits_{|x-e|+|x|\leq 2\tau}\frac{q(x)}{|x||x-e|}dx&=\frac{1}{2}\int\limits_{1}^{2\tau}\int\limits_{0}^{\pi}\int\limits_{0}^{2\pi}q(\ln(r+\sqrt{r^{2}-1}),\theta,\phi)\sin\phi d\theta d\phi dr.
	\end{align*}
	
	Now using $\eqref{expression for Q}$, we get
	\begin{equation}\label{estimate for volume integral in prolate}
	\begin{aligned}
	\int\limits_{|x-e|+|x|\leq 2\tau}\frac{q(x)}{|x||x-e|}dx\leq C \int\limits_{1}^{2\tau}Q(r)dr. 
	\end{aligned}
	\end{equation}
	Now applying this inequality in Equation \eqref{Estimate} and noting that $q(x)=q_{1}(x)-q_{2}(x)\geq 0$, we {\color{black}have}
	\begin{equation}\label{estimate in prolate}
	\begin{aligned}
	Q(2\tau)\leq CK\int\limits_{1}^{2\tau}Q(r)dr.
	\end{aligned}
	\end{equation}
	Now Equation $\eqref{estimate in prolate}$ holds for all $\tau\in [1/2 ,T/2]$ and since $q(x)\geq 0$, for all $x\in\rr^{3}$, {\color{black}therefore using the Gronwall's inequality,  we get}
	\begin{align*}
	Q(2\tau)=0,\ \ \ \tau\in[1/2,T/2].
	\end{align*}
	Now from Equation \eqref{Definition of Q}, again using $q(x)\geq 0$,  we have $q(x)=0$, for all $x\in\rr^{3}$ such that $|x|+|x-e|\leq T$. The proof is complete.
	\subsection{Proof of Theorem \ref{Separated Result 2}}
	Again first, we consider the surface integral in $\eqref{Estimate}$ and denote it by $Q(2\tau)$:
	\begin{align*}
	Q(2\tau):=\int\limits_{|x|+|x-e|=2\tau}\frac{q(x)}{|2\tau x-|x|e|}dS_{x}.
	\end{align*}
	and $q(x):=a(|x|+|x-e|)$.
	Now from Equations $\eqref{Prolate}$, $\eqref{Surface and Volume measure on ellipsoid}$ and \eqref{expression for Q} and hypothesis $q_{i}(x)=a_{i}(|x|+|x-e|)$ of the theorem, we get
	\begin{equation}\label{final expression for Q}
	\begin{aligned}
	Q(2\tau)=\frac{1}{8}\int\limits_{0}^{\pi}\int\limits_{0}^{2\pi}a(2\tau)\sin\phi d\theta d\phi=\frac{\pi}{2}a(2\tau).
	\end{aligned}
	\end{equation}
	Now consider the integral 
	\begin{align*}
	\int\limits_{|x|+|x-e|\leq 2\tau}\frac{q(x)}{|x||x-e|}dx.
	\end{align*}
	Again using $\eqref{Prolate}$ and $\eqref{Surface and Volume measure on ellipsoid}$ in the above integral, we have 
	\begin{align*}
	\begin{aligned}
	\int\limits_{|x-e|+|x|\leq 2\tau}\frac{q(x)}{|x||x-e|}dx&= \frac{1}{2}\int\limits_{\cosh\rho\leq 2\tau}\int\limits_{0}^{\pi}\int\limits_{0}^{2\pi}q(\rho,\theta,\phi)\sinh\rho\sin\phi d\theta d\phi d\rho.\\
	\end{aligned}
	\end{align*} 
	After substituting $\cosh\rho=r$ and $\rho=\ln(r+\sqrt{r^{2}-1})$,  we get 
	\begin{align*}
	\left\lvert \ \int\limits_{|x-e|+|x|\leq 2\tau}\frac{q(x)}{|x||x-e|}dx\right\rvert&=\left\lvert\frac{1}{2}\int\limits_{1}^{2\tau}\int\limits_{0}^{\pi}\int\limits_{0}^{2\pi}a(r)\sin\phi d\theta d\phi dr\right\rvert\\
	&\leq C\int\limits_{1}^{2\tau}|a(r)|dr.
	\end{align*}
	
	Now using this inequality and Equation $\eqref{final expression for Q}$ in $\eqref{Estimate}$,  we see
	\begin{equation}\label{estimate in Prolate for ellipsoidal}
	|a(2\tau)|\leq C\int\limits_{1}^{2\tau}|a(r)|dr.
	\end{equation}
	Now Equation $\eqref{estimate in Prolate for ellipsoidal}$ holds for all $\tau\in[1/2,T/2]$, so using the Gronwall's inequality, we have 
	\begin{align*}
	a(2\tau)=0,\  \mbox{for}\ \tau\in[1/2,T/2].
	\end{align*}
	Thus,
	we have $q(x)=0$, for all $x\in\rr^{3}$ such that $|x|+|x-e|\leq T$. {\color{black} This conclude the proof of Theorem \ref{Separated Result 2}.}
	{\section{Conclusion}
		In this paper, we studied an inverse problem for the wave equation with \emph{separated} point source and point receiver data. Our approach is based on construction  of spherical wave solution using the solution to a Goursat problem, combined with the solution to an adjoint problem, we ended up with an integral identity. Then using the prolate-spheroidal co-ordinates and Grownwall's inequality, we completed the proof of the main results.}
	
	\section*{Acknowledgement}
	{\color{black}The author thanks the anonymous referees for useful comments which helped him to improve the paper.} The author would like to thank his advisor Venky Krishnan for his great motivation and useful discussions. He would like to thank Prof. Rakesh for suggesting this problem during the workshop ``Advanced Instructional School on Theoretical and Numerical Aspects of Inverse Problems, June 16--28, 2014'' held at TIFR Centre for Applicable Mathematics, Bangalore, India, and for suggesting the use of solution to the adjoint problem. He also would like to thank Prof. Paul Sacks for stimulating discussions.


\begin{thebibliography}{99}
		
		
		\vspace*{1mm}
		\footnotesize
		
		
		\bibitem{Bube and Burridge 1-D} K. P. Bube and R. Burridge;
		The one-dimensional inverse problem of reflection seismology.
		SIAM Rev. 25 (1983), no. 4, 497--559.
		\bibitem{Bukhgem and Klibanov uniqueness} A. L. Bukhge\~{i}m and M.V. Klibanov; Uniqueness in the large of a class of multidimensional inverse problems. (Russian)
		Dokl. Akad. Nauk SSSR 260 (1981), no. 2, 269--272.
		\bibitem{Burridge integral equation} R. Burridge; The Gel'fand-Levitan, the Marchenko, and the Gopinath-Sondhi integral equations of inverse scattering theory, regarded in the context of inverse impulse-response problems. Wave Motion 2 (1980), no. 4, 305--323. 
		\bibitem{Friedlander book} F.G.  Friedlander; The wave equation on a curved space-time. Cambridge University Press, Cambridge-New York-Melbourne, 1975. Cambridge Monographs on Mathematical Physics, No. 2.
		
		
		\bibitem{Klibanov partial point source} M.V. Klibanov; Some inverse problems with a `partial' point source. Inverse Problems 21 (2005), no. 4, 1379--1390.
		
		\bibitem{LDS} M.M. Lavrent\'ev;  E.Yu.  Derevtsov and V.A. Sharafutdinov; "On the determination of an
		optical body, situated in a homogeneous medium, from its images," Dokl. Akad. Nauk SSSR,
		260, No. 4, 799-803 (1981).
		
		\bibitem{Levrentev Romanov and Shishatski book ill-posed} M.M. Lavrent'ev; V.G. Romanov and S.P. Shishat$\cdot$ski\~{i}; Ill-posed problems of mathematical physics and analysis. Translated from the Russian by J. R. Schulenberger. 
		\bibitem{Li partial point source} S. Li; Estimation of coefficients in a hyperbolic equation with impulsive inputs. J. Inverse Ill-Posed Probl. 14 (2006), no. 9, 891--904. 
			\bibitem{Rakesh impedance inversion} Rakesh; An inverse impedance transmission problem for the wave equation. Comm. Partial Differential Equations 18 (1993), no. 3-4, 583--600. 
				\bibitem{Rakesh inversion of spherically symmetric potential} Rakesh; Inversion of spherically symmetric potentials from boundary data for the wave equation. Inverse Problems 14 (1998), no. 4, 999--1007.
			\bibitem{Rakesh Webster's equation} Rakesh; Characterization of transmission data for Webster's horn equation. Inverse Problems 16 (2000), no. 2, L9--L24. 
			\bibitem{Rakesh single-coincident source-receiver pair} Rakesh; Inverse problems for the wave equation with a single coincident source-receiver pair. Inverse Problems 24 (2008), no. 1, 015012, 16 pp.
			\bibitem{Rakesh layered medium point source}Rakesh; An inverse problem for a layered medium with a point source.  Problems 19 (2003), no. 3, 497--506.	
			\bibitem{Rakesh and Sacks impedance inversion} Rakesh and P. Sacks; Impedance inversion from transmission data for the wave equation. Wave Motion 24 (1996), no. 3, 263--274. 
			
			\bibitem{Rakesh and Sacks uniqueness angular controlled potential}Rakesh and P. Sacks;  Uniqueness for a hyperbolic inverse problem with angular control on the coefficients. J. Inverse Ill-Posed Probl. 19 (2011), no. 1, 107--126.
			
		\bibitem{Rakesh and Symes uniqueness result} Rakesh; W.W. Symes; Uniqueness for an inverse problem for the wave equation. Comm. Partial Differential Equations 13 (1988), no. 1, 87--96.	
	
		 
		\bibitem{Rakesh and Uhlmann inverse backscattering}Rakesh and G. Uhlmann; The point source inverse back-scattering problem. Analysis, complex geometry, and mathematical physics: in honor of Duong H. Phong, 279--289. 
		\bibitem{Romanov}V.G. Romanov; On the problem of determining the coefficients in the lowest order terms of a hyperbolic equation. (Russian. Russian summary) Sibirsk. Mat. Zh. 33 (1992), no. 3, 156--160, 220; translation in Siberian Math. J. 33 (1992), no. 3, 497--500. 
		\bibitem{Romanov book integral geometry}V.G. Romanov; Integral geometry and inverse problems for hyperbolic
		equations, volume 26. Springer Science and Business Media, 2013.
		\bibitem{Romanov problem of determining the two coefficients} V.G. Romanov and D.I.  Glushkova; The problem of determining two coefficients of a hyperbolic equation. (Russian) Dokl. Akad. Nauk 390 (2003), no. 4, 452--456. 
			
		\bibitem{Sacks and Symes uniqueness and continuous dependence}  P. Sacks and W.W. Symes; Uniqueness and continuous dependence for a multidimensional hyperbolic inverse problem. Comm. Partial Differential Equations 10 (1985), no. 6, 635--676. 
		\bibitem{Santosa Syemes} F. Santosa and W.W. Symes;  High-frequency perturbational analysis of the surface point-source response of a layered fluid. J. Comput. Phys. 74 (1988), no. 2, 318--381.		
		\bibitem{Sondhi} M.M. Sondhi; A survey of the vocal tract inverse problem: theory, computations and experiments. Inverse problems of acoustic and elastic waves (Ithaca, N.Y., 1984), 1--19, SIAM, Philadelphia, PA, 1984.
		
		
		
		
		
		\bibitem{Stefanov uniqueness} P.D. Stefanov; A uniqueness result for the inverse back-scattering problem. Inverse Problems 6 (1990), no. 6, 1055--1064.
		
		
		
		
		\bibitem{Symes} W.W. Symes; The seismic reflection inverse problem. Inverse Problems 25 (2009), no. 12, 123008, 39 pp.
		
		
		
		
	\end{thebibliography}
\end{document}